\documentclass[letterpaper, 11pt]{article}

\usepackage{graphics,graphicx,xcolor,float,enumerate, bm}
\usepackage{multicol} 
\usepackage{parskip}
\usepackage{amsmath,amssymb}
\usepackage{multirow}
\usepackage[utf8]{inputenc}
\usepackage{fancyhdr}
\usepackage[title]{appendix}
\usepackage{wasysym}
\usepackage{url}
\usepackage{mathtools}
\usepackage{booktabs}
\usepackage{float}

\usepackage[font=footnotesize,labelfont=small]{caption}
\RequirePackage{geometry}
\newtheorem{remark}{Remark}

\newcommand{\vct}[1]{\bm{#1}}
\newcommand{\mtx}[1]{\mathsf{#1}}

\newcommand{\vtwo}[2]{\left[\begin{array}{c} #1 \\ #2 \end{array}\right]}
\newcommand{\vthree}[3]{\left[\begin{array}{c} #1 \\ #2 \\ #3 \end{array}\right]}

\newcommand{\mtwo}[4]{\left[\begin{array}{cc}    #1 & #2 \\ #3 & #4  \end{array}\right]}
\newcommand{\mthree}[9]{\left[\begin{array}{ccc} #1 & #2 & #3 \\
                                                 #4 & #5 & #6 \\
                                                 #7 & #8 & #9 \end{array}\right]}

  \newtheorem{definition}{Definition}

\numberwithin{pro}{section}

\title{A fast direct solver for two dimensional quasi-periodic multilayered  media scattering problems, Part II}
\author{Yabin Zhang\thanks{Department of Mathematics, University of Michigan, Ann Arbor} 
\;and Adrianna Gillman\thanks{Department of Applied Mathematics, University of Colorado, Boulder}}

\date{}

\begin{document}
\maketitle

\begin{abstract}
This manuscript is the second in a series presenting fast direct solution techniques for solving 
two-dimensional wave scattering problems 
from quasi-periodic multilayered structures. 
The fast direct solvers presented in the series are for the linear system that results from the  discretization of a robust integral formulation.  The fast direct solver presented in this manuscript has a computational cost that scales linearly with respect to the number of discretization points on the interfaces and the number of 
layers.  The latter is an improvement over the previous solver and makes the new solver more efficient especially for problems involving multiple incident angles and changes to the layered media.  Numerical results illustrate the improved performance of the new solver over the previous one.
\end{abstract}

\section{Introduction}
\label{sec:intro}
\begin{figure}[htbp]
\centering
\begin{picture}(200,250)(100,0)
\put(0,0){\includegraphics[trim={0 2cm 0 2cm},clip,scale=0.7]{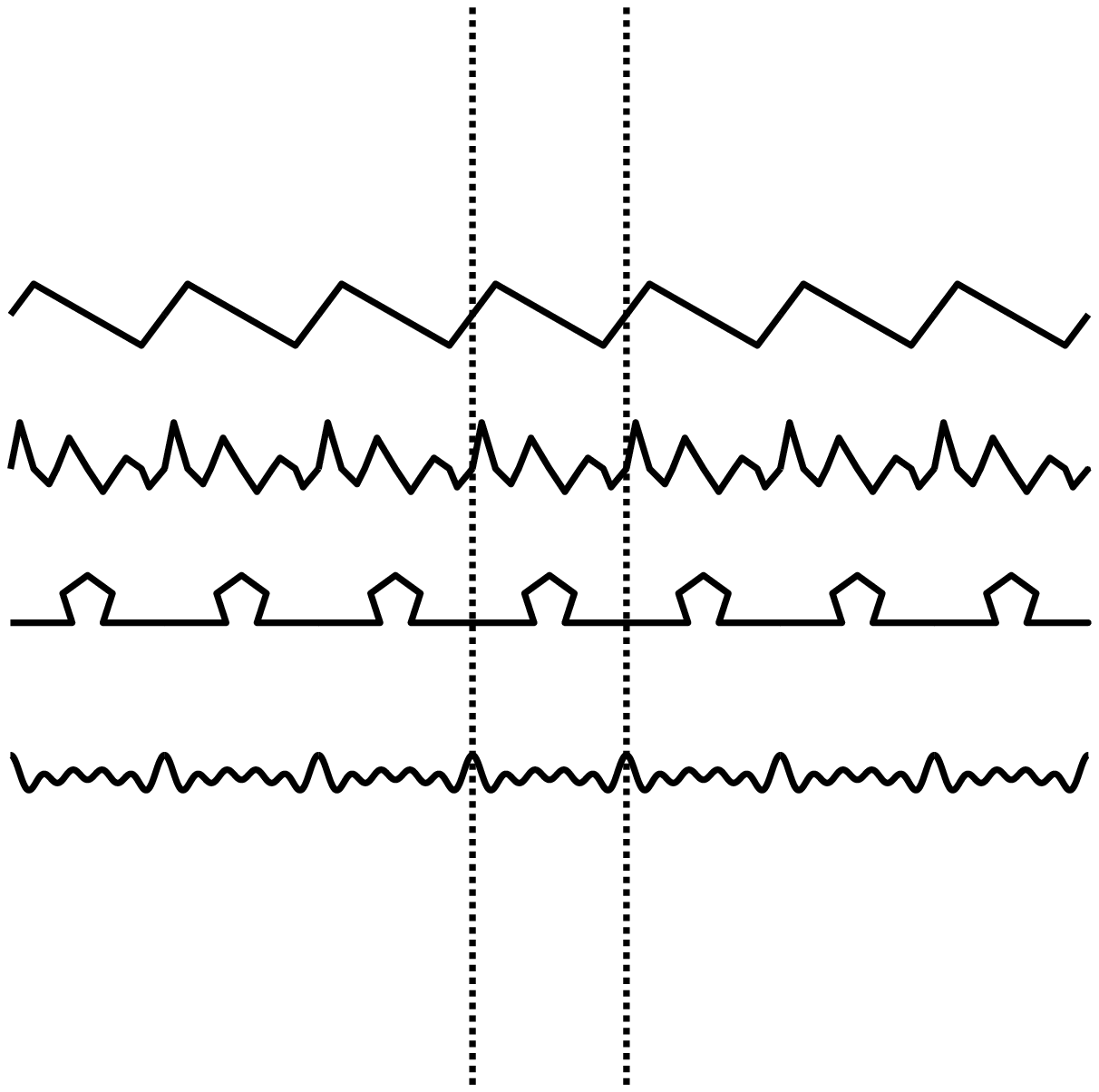}}
\put(100,220){\vector(1,-1){20}}
\put(100,210){\line(1,1){10}}
\put(105,205){\line(1,1){10}}
\put(110,230){$u^{inc}$}
\put(280,190){$\Omega_1$}
\put(280,145){$\Omega_2$}
\put(280,110){$\Omega_3$}
\put(280,75){$\Omega_4$}
\put(280,35){$\Omega_5$}
\put(330,160){$\Gamma_1$}
\put(330,130){$\Gamma_2$}
\put(330,90){$\Gamma_3$}
\put(330,60){$\Gamma_4$}
\put(196,200){\vector(-1,0){10}}
\put(210,200){\vector(1,0){10}}
\put(201,198){$d$}
\end{picture}
\begin{center}
\caption{A 5-layer(4-interface) periodic geometry. 7 periods are shown.}
\label{fig:five_layer_geom}
\end{center}
\end{figure}

This manuscript considers the quasi-periodic scattering problem defined in layered media by the following Helmholtz problem: \begin{equation}
 \label{eq:basic}
 \begin{split}
(\Delta + \omega_i^2) u_i(\vct{x}) &= 0  \  \qquad \vct{x}\in \Omega_i\\
 u_1-u_2 &= -u^{\rm inc}(\vct{x}) \  \qquad \vct{x}\in \Gamma_1\\
 \frac{\partial u_1}{\partial\nu} -\frac{\partial u_2}{\partial \nu} &= -\frac{\partial u^{\rm inc}}{\partial \nu} \  \qquad \vct{x}\in \Gamma_1\\
u_{i}-u_{i+1} &= 0 \  \qquad \vct{x}\in \Gamma_{i}, \ 1<i<I+1\\
 \frac{\partial u_i}{\partial\nu} -\frac{\partial u_{i+1}}{\partial \nu} &= 0 \ \qquad  \vct{x}\in \Gamma_{i},\ 1<i<I+1
 \end{split}
\end{equation}
\noindent
where $u_i$ is the unknown solution in the region $\Omega_i\in\mathbb{R}^2$,
 the wave number in $\Omega_i$ is given by $\omega_i$ for $i = 1,\ldots,I+1$,
 and $\nu(\vct{x})$ is the normal vector at $\vct{x}$.
The interface $\Gamma_i$ for $i =1,\ldots,I$ between each layer is periodic
with period $d$. 
Figure \ref{fig:five_layer_geom} illustrates an example of a 5-layer periodic geometry.

The multilayered media problem (\ref{eq:basic}) can be transformed into a collection 
of integral equations defined on each interface via the formulation in \cite{multilayer_part1,Cho:15}.
This formulation is robust even at Wood's anomalies and is amenable to fast direct solvers.  The basic
idea is to separate one period of the geometry into two pieces: a unit cell (box) containing the layers and the infinite vertical strips outside of that box.  Inside the unit cell, the solution is expressed as an integral operator defined on a single period plus integral operators on 
neighboring periods and an additional term for enforcing quasi-periodicity. In the vertical strips 
outside of the unit cell, the solution is represented via a Rayleigh-Bloch expansion.  
The fast direct solver from \cite{multilayer_part1} utilizes a classic fast direct solver
for integral equations  \cite{2012_martinsson_FDS_survey,2012_ho_greengard_fastdirect,2013_3DBIE} for the operator defined on a unit cell 
and low-rank factorizations for the interactions with the neighboring periodic copies
as in \cite{2013_Gillman,2016_periodic_stokes}.  The resulting solver scales linearly
with the number of discretization points on the interfaces and is able to efficiently 
handle problems involving a small-to-moderate number of complicated interfaces on 
a desktop computer.  It is also able to efficiently handle multiple incident angles.
Unfortunately, the solver from \cite{multilayer_part1} is not optimal for problems involving
local changes in the 
layered medium such as changes in wave number or interface geometries.  Such problems arise
in applications such as optimal design problems.
This manuscript presents a new direct solver which is optimal for these problems making the integral equation formulation from \cite{Cho:15} practical for problems involving many local changes in the medium.  



\subsection{High level view of the solution technique}
The linear system that results from the integral equation formulation in \cite{Cho:15} results in a 
block linear system where each block row equation enforces a different part of the requirements
of the integral formulation: boundary conditions through the interfaces, the periodicity and the 
radiation condition.  Following the solution technique in \cite{Cho:15}, the solution to the 
block system can be constructed in block form requiring the inverse of a Schur complement 
operator that is block tri-diagonal 
where the off-diagonal blocks are low-rank and the 
diagonal blocks can be written as the sum of a full-rank matrix with low-rank matrices.  Thus the inverse can be applied via a Woodbury formula.  Each 
full rank diagonal block matrix corresponds to the discretization of an integral operator on each 
interface.  These matrices can be inverted quickly via a fast direct solver such as the \textit{Hierarchically
Block Separable (HBS)} \cite{2012_martinsson_FDS_survey,2012_ho_greengard_fastdirect,2013_3DBIE},
the \textit{Hierarchically Semi-Separable (HSS)}
 \cite{2009_xia_superfast,2007_shiv_sheng,2010_xia}, the \textit{Hierarchical interpolative
 factorization (HIF)} \cite{2014_HIF}, and the $\mathcal{H}$ and $\mathcal{H}^{2}$-matrix methods
 \cite{2010_borm_book,2004_borm_hackbusch}.  The low-rank matrices in the Schur complement are made up of 
 a sum of low-rank matrices corresponding to the interaction of an interface with its neighbors 
 and matrices that help enforce the periodicity and radiation conditions.  
The fast solver in this manuscript constructs the low-rank factors in the same manner as 
in \cite{multilayer_part1}.  

The fast inverse of the discretized boundary integral equations on the diagonal blocks allows the solver to scale linearly with respect to the number of unknowns placed on each interface, and the Schur complement formulation allows the solver to scale linearly with respect to the total number of interfaces. The fast direct solver
in this manuscript is ideal for optimal design problems,
since the additional cost for changing an interface geometry and/or a layer wave number scales linearly with the number of discretization points on the affected interfaces instead of the total number of unknowns on the entire structure. Once the solver is built, new incident angles can also be handled with small extra cost as the most expensive calculations in building the solver is independent of incident angle or Bloch phase.

\subsection{Outline}
The manuscript begins by briefly reviewing the linear system that results from the discretization of the
robust integral formulation for equation (\ref{eq:basic}) from \cite{Cho:15} in Section \ref{sec:per}.
 Next the new fast direct solver is presented in Section \ref{sec:FDS}.  Section \ref{sec:extensions} describes extensions of the proposed solution technique that make it useful for practical applications.  
Section \ref{sec:numerics} illustrates the performance of the solver for a selection of test problems.
Finally, Section \ref{sec:summary} concludes the manuscript and discusses future directions.

\section{The linear system and Schur complement}
\label{sec:per}
\begin{figure}[h]
\centering
\begin{picture}(200,260)(100,-10)
\put(-50,0){\includegraphics[trim={0 0 0 0},clip,scale=0.5]{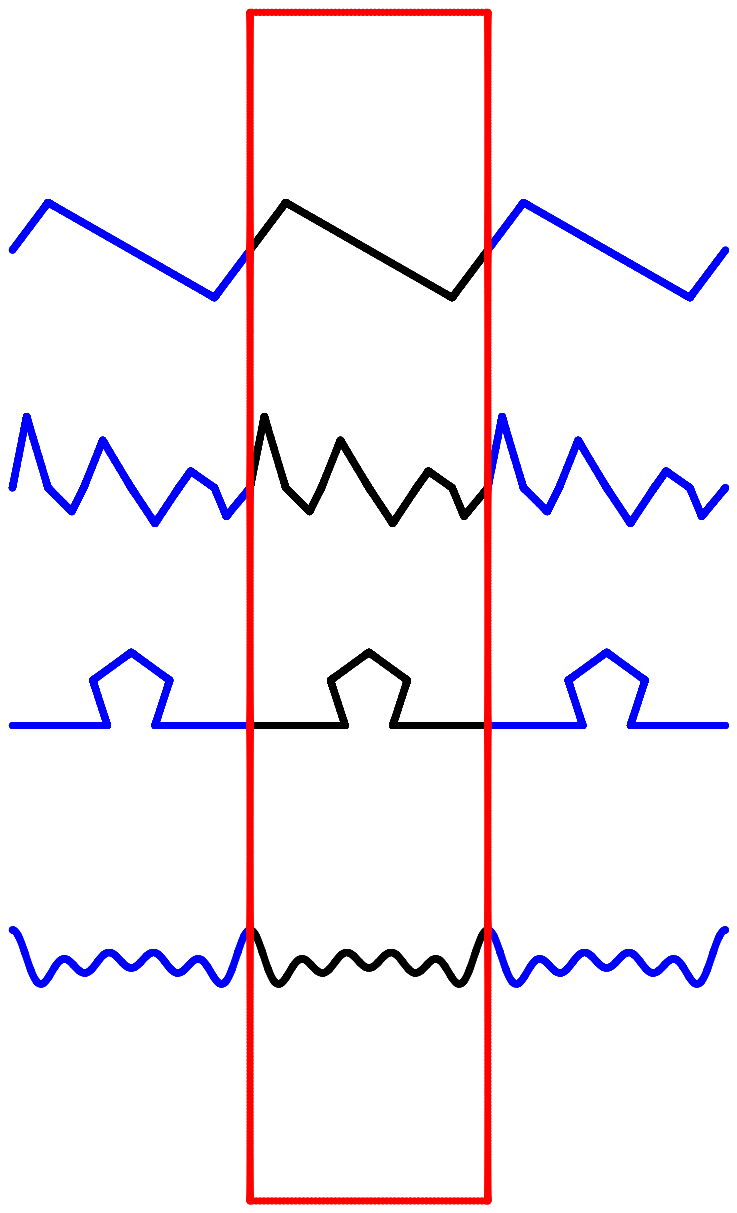}}
\put(90,-10){(a)}
\put(67,180){\color{red} $L$}
\put(118,180){\color{red} $R$}
\put(80,200){\color{red} $y=y_U$}
\put(80,10){\color{red} $y =y_D$}
\put(90,166){\footnotesize$\Gamma_1$}
\put(90,135){\footnotesize$\Gamma_2$}
\put(90,105){\footnotesize$\Gamma_3$}
\put(90,63){\footnotesize$\Gamma_4$}
\put(180,-30){\includegraphics[trim={5cm 0 5cm 0},clip,scale=0.65]{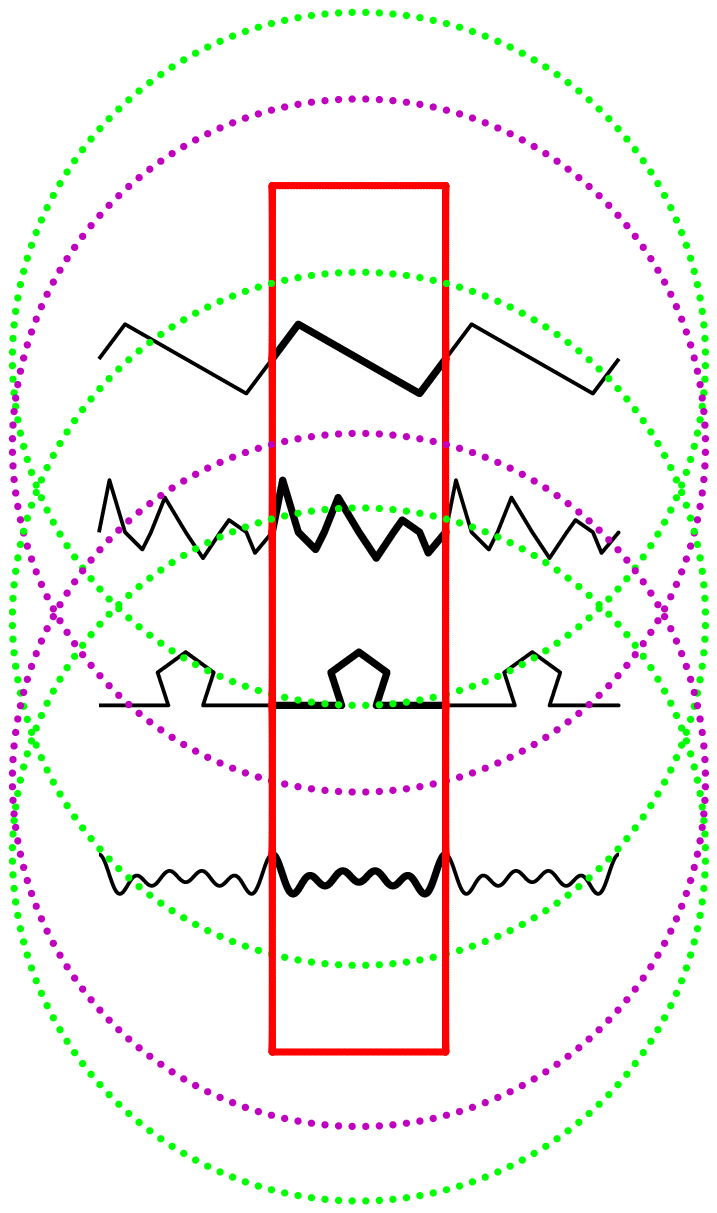}}
\put(280,-10){(b)}
\put(242,170){\footnotesize$\Omega_1$}
\put(242,140){\footnotesize$\Omega_2$}
\put(242,110){\footnotesize$\Omega_3$}
\put(242,75){\footnotesize$\Omega_4$}
\put(242,38){\footnotesize$\Omega_5$}
\put(300,222){\color{green}$P_1$}
\put(310,182){\color{magenta}$P_2$}
\put(345,112){\color{green}$P_3$}
\put(310,35){\color{magenta}$P_4$}
\put(303,-2){\color{green}$P_5$}

\end{picture}
\begin{center}
\caption{
This figure illustrates a 5-layer periodic geometry with artificial walls and proxy circles. 
Only three periods of the infinite periodic geometry are shown.  The period contained within the
unit cell is plotted in black while neighboring periods are plotted in blue.
Figure (a) illustrates the notation for the unit cell with left, right, upper, and lower boundary $L$, $R$, $y_U$, and $y_D$ shown in red lines.
Figure (b) illustrates the proxy circles $P_i$ for each layer.  The color of the proxy circles alternates 
between green and magenta.}
\label{fig:five_layer_geom_unitcell}
\end{center}
\end{figure}

The integral formulation from \cite{Cho:15} avoids using the quasi-periodic Green's function by introducing auxiliary unknowns on one period of the layered structure, referred to as the unit cell, to 
enforce periodicity in the horizontal direction. 
Figure \ref{fig:five_layer_geom_unitcell}(a) illustrates an example of a unit cell.  
Inside the unit cell, the solution in each layer is represented via an integral equation defined on the interfaces plus a collection of point charges 
of unknown magnitude to capture the quasi-periodicity of the approximate solution.  These point charges are placed on a proxy circle 
that encloses the layer in the unit cell as illustrated in Figure \ref{fig:five_layer_geom_unitcell}(b).  
Outside the unit cell (in the positive and negative $y-$direction) the solution is represented via a Rayleigh-Bloch expansion which naturally satisfies the radiation condition. So there are three sets of unknowns: the boundary charge densities on each interface, the magnitude of the 
point charges for capturing the quasi-periodicity of the solution and the coefficients of the Rayleigh-Bloch expansion.  These unknowns
are found by enforcing continuity of the solution and the flux through the interfaces, enforcing the quasi-periodicity of the solution 
on the left and right walls of the unit cell 
and 
enforcing the continuity of the two solution representations and the flux 
 through the top and bottom of the unit cell.
The left and right walls of the unit cell are labeled $L$ and $R$, respectively, in Figure \ref{fig:five_layer_geom_unitcell}(a).  
The top and bottom walls of the unit cell are labeled $y=y_U$ and $y=y_D$, respectively, in Figure \ref{fig:five_layer_geom_unitcell}(a).
Detailed descriptions of the integral formulation are presented in \cite{Cho:15,multilayer_part1}. 


The linear system that results from the discretization of the integral formulation in \cite{Cho:15} 
has the following rectangular form:
\begin{equation}
 \mthree{\mtx{A}}{\mtx{B}}{\vct{0}}{\mtx{C}}{\mtx{Q}}{\vct{0}}{\mtx{Z}}{\mtx{V}}{\mtx{W}}\vthree{\vct{\hat{\sigma}}}{\vct{c}}{\vct{a}}
 = \vthree{\vct{f}}{\vct{0}}{\vct{0}}
 \label{eq:bigblock}
\end{equation}
where 
\begin{itemize}
    \item $\mtx{A}$ is a matrix of size $2N\times 2N$ that results from the discretization of integral operators on the interfaces where $N=\sum_{i=1}^I N_i$ is the total number of discretization points on all the interfaces and $N_i$ is the number of discretization points on interface $i$,
\item
$\mtx{B}$ is a matrix of size $2N\times P$ where $P = \sum_{l=1}^{I+1} P_l$ is the total number of point charges placed on the proxy circles for all layers and $P_l$ is the number of number of point charges for layer $l$, 
\item
$\mtx{C}$ is a matrix of size $2(I+1) M_w\times 2N$ and
$\mtx{Q}$ is a matrix of size $2(I+1)M_w\times P$ where $M_w$ is the number of locations per layer on the left 
and right vertical walls of the unit cell where quasi-periodicity is enforced,
\item 
$\mtx{Z}$ is a matrix of size $4M\times 2N$, 
$\mtx{V}$ is a matrix of size $4M\times P$,
and $\mtx{W}$ is a matrix of size $4M\times 2(2K+1)$ where $K$ is the number of terms used in the Rayleigh-Bloch 
expansion and $M$ is the number of locations on the top and bottom of the unit cell where the continuity of the solution
is enforced with outside the unit cell.
\end{itemize}
The first row equation in (\ref{eq:bigblock}) enforces continuity of the solution through the interfaces.
The second row equation in (\ref{eq:bigblock}) enforces quasi-periodicity of the solution.  The 
last row equation enforces continuity of the solution outside the unit cell.  Detailed definitions of the matrix block entries are given in Section  2.3 of \cite{multilayer_part1}.

The unknowns that need to be found are:  the boundary charge densities $\hat{\bm{\sigma}}$ on all of the interfaces, 
the coefficients of the point charges (basis functions) $\bm{c}$ defined on the (uniformly) sampled locations of the proxy circles,
 and the coefficients of the Rayleigh-Block expansions $\bm{a}$. 
The block $\bm{f}$ of the right-hand-side vector contains zero in all the entries except for the ones corresponding
to the top interface whose value is determined by the incident plane wave of the top layer.

The matrix blocks in (\ref{eq:bigblock}) are all block sparse. A reordering of the unknowns allows for the
solution of this linear system to be written as a block solve involving a collection of block tri-diagonal matrices.
Specifically, we introduce the vector $\vct{x}$ which is the ordering of the vector $\begin{bmatrix}\vct{c}\\ \vct{a} \end{bmatrix}$ defined by 
\begin{equation}
\label{eq:reorder}
\vct{x}=\begin{bmatrix}
\vct{x}_1\\
\vdots\\
\vct{x}_{I+1}
\end{bmatrix}
\mbox{ with }
\vct{x}_1=\begin{bmatrix}
\vct{c}_1\\
\vct{a}^U
\end{bmatrix},
\vct{x}_i=\vct{c}_i, \mbox{ for }i=2,\dots, I, 
\mbox{ and }
\vct{x}_{I+1}=\begin{bmatrix}
\vct{c}_{I+1}\\
\vct{a}^D
\end{bmatrix}.
\end{equation} where 
the subscript indicates which layer the unknowns on the proxy circle belong to  (e.g, $\bm{c}_1$ are the unknowns placed on the proxy circle of the top layer) and the vectors $\bm{a}^U$ and $\bm{a}^D$ denote the unknown coefficients of the Rayleigh-Bloch expansions defined on ${y}=y_U$ and $y=y_D$ respectively.  Then the block linear system (\ref{eq:bigblock}) is written as 
\begin{equation}
 \mtwo{\mtx{A}}{\hat{\mtx{B}}}{\hat{\mtx{C}}}{\hat{\mtx{Q}}}\vtwo{\vct{\hat{\sigma}}}{\vct{x}}
 = \vtwo{\vct{f}}{\vct{0}}.
 \label{eq:rearrangedblock}
\end{equation}
where 
$\hat{\mtx{B}}$ is the reordered version of $\mtx{B}$, $\hat{\mtx{C}}$ is the reorder version of 
$\begin{bmatrix}
\mtx{C}\\
\mtx{Z}
\end{bmatrix}$ and $\hat{\mtx{Q}}$ is the reordered version of $ \begin{bmatrix}
\mtx{Q} & \mtx{0}\\
\mtx{V} & \mtx{W}
\end{bmatrix}$.  

We choose to solve this system in a block solve format;
\begin{equation}
    \begin{array}{rl}
    \vct{x} & = -\hat{\mtx{Q}}^{\dagger} \hat{\mtx{C}}\hat{\vct{\sigma}}\\
    \hat{\vct{\sigma}} & = \mtx{S}^{-1}\vct{f}
    \end{array}
    \label{eq:blocksolve}
\end{equation}
where 
\begin{equation}\label{eq:schur1}
\mtx{S} = \mtx{A}-\hat{\mtx{B}}\hat{\mtx{Q}}^{\dagger} \hat{\mtx{C}}
\end{equation}
denotes
the Schur complement of the block system.  Thanks to the 
reordering of the unknowns the matrices $\mtx{S}$ and $\hat{\mtx{Q}}$ are block tri-diagonal.  This
allows for the processing of the solve to be very fast.  The block solve in (\ref{eq:blocksolve}) is the same as in \cite{Cho:15}.  

The remainder of this section provides a high level view of the different matrices that make up the linear system 
(\ref{eq:rearrangedblock}).
Section \ref{sec:blockA} reports the 
tri-diagonal block entries of the matrix $\mtx{A}$.  Section \ref{sec:otherblocks} reports 
 the entries of the other block matrices in the reordered linear system.



\subsection{Block structure of $\mtx{A}$}
\label{sec:blockA}
The matrix $\mtx{A}$ in equation (\ref{eq:bigblock}) is (almost always) the largest matrix in the linear system.  Fortunately $\mtx{A}$ has structure that makes it amenable to accelerated linear algebra.  We first note that $\mtx{A}$ is a block tri-diagonal matrix
$$
{\mtx{A}}=
\begin{bmatrix}
{\mtx{A}}_{1,1} & \mtx{A}_{1,2} & 0 & \cdots &0 &  0 &  0 \\
\mtx{A}_{2,1} &  \mtx{A}_{2,2} &  \mtx{A}_{2,3} & \cdots  &  0&0& 0\\
0 & \mtx{A}_{3,2} &\mtx{A}_{3,3} & \cdots  &  0&0& 0\\
\vdots & \vdots &\vdots &\vdots &\vdots & \vdots & \vdots\\
0 & 0 & 0 & \cdots & \mtx{A}_{I-1,I-2}& \mtx{A}_{I-1,I-1} & \mtx{A}_{I-1,I}\\
0 & 0 & 0 & \cdots & 0& \mtx{A}_{I,I-1} &{\mtx{A}}_{I,I}\\
\end{bmatrix} 
$$
where the diagonal block $\mtx{A}_{i,i}$ corresponds to the self interactions of the $i^{\rm th}$ 
interface in the unit cell and its interaction with the adjacent left and right neighbor periods and the off-diagonal block $\mtx{A}_{i,j}$ corresponds to the interaction between the $i^{\rm th}$ interface and the interface above when $j=i-1$ or the interface below when $j=i+1$.  For simplicity of presentation, we denote the diagonal blocks as follows:
$$\mtx{A}_{i,i} = \mtx{A}_{i,i}^s + \mtx{A}_{i,i}^{pm}$$
where $\mtx{A}_{i,i}^s$ denotes the interaction of interface $i$ inside the unit cell with itself and $\mtx{A}_{i,i}^{pm}$ denotes the interaction between 
the $i^{\rm th}$ interface in the unit cell and the left and right copies of it.  The superscripts ``$s$'' and ``$pm$'' stand for ``self" and  the left and right period (plus and minus) interactions respectively. This matches notation in \cite{multilayer_part1}.

With this notation, the matrix $\mtx{A}$ can be written as
\begin{equation}
{\mtx{A}}=
\mtx{A}_0+
%
\begin{bmatrix}
{\mtx{A}}^{pm}_{1,1} & \mtx{A}_{1,2} & 0 & \cdots &0 &  0 &  0 \\
\mtx{A}_{2,1} &  \mtx{A}^{pm}_{2,2} &  \mtx{A}_{2,3} & \cdots  &  0&0& 0\\
0 & \mtx{A}_{3,2} &\mtx{A}^{pm}_{3,3} & \cdots  &  0&0& 0\\
\vdots & \vdots &\vdots &\vdots &\vdots & \vdots & \vdots\\
0 & 0 & 0 & \cdots & \mtx{A}_{I-1,I-2}& \mtx{A}^{pm}_{I-1,I-1} & \mtx{A}_{I-1,I}\\
0 & 0 & 0 & \cdots & 0& \mtx{A}_{I,I-1} &{\mtx{A}}^{pm}_{I,I}\\
\end{bmatrix} 
\label{eq:def_Aplus}
\end{equation}
where $\mtx{A}_0$ is a block diagonal matrix
\begin{equation}
\label{equ:def_Anot}
\mtx{A}_0=
\begin{bmatrix}
{\mtx{A}}^s_{1,1} & 0 & 0 & \cdots &0 &  0 &  0 \\
0 &  \mtx{A}^s_{2,2} & 0 & \cdots  &  0&0& 0\\
0 & 0 &\mtx{A}^s_{3,3} & \cdots  &  0&0& 0\\
\vdots & \vdots &\vdots &\vdots &\vdots & \vdots & \vdots\\
0 & 0 & 0 & \cdots & 0& \mtx{A}^s_{I-1,I-1} &0\\
0 & 0 & 0 & \cdots & 0& 0 &{\mtx{A}}^s_{I,I}\\
\end{bmatrix}.
\end{equation}
All the non-self interaction matrices, including interactions between the neighboring periods of an interface and the interactions between vertically neighboring interfaces, are low-rank. For example 
$\mtx{A}_{i,i}^{pm}$ and $\mtx{A}_{1,2}$ are low-rank.
This means that each of the blocks in the block tri-diagonal matrix in the right hand side of equation (\ref{eq:def_Aplus}) are low-rank.

\subsection{The block structure of the remaining matrices in the reordered linear system}
\label{sec:otherblocks}
The rearrangement of the auxiliary unknowns $\vct{c}$ and $\vct{a}$ via the 
ordering defined in equation (\ref{eq:reorder}) changes
the sparsity pattern of the non-principal block matrices
in (\ref{eq:rearrangedblock}).  This section presents the 
sparsity and entries of the reordered matrices.

The matrix $\hat{\mtx{B}}$  is $\begin{bmatrix} \mtx{B} & \mtx{0} \end{bmatrix}$ with its columns reordered according to the new ordering of the unknowns and is defined as follows
$$\hat{\mtx{B}}=
\begin{bmatrix}
\hat{\mtx{B}}_{1,1} & \mtx{B}_{1,2} & 0 & \cdots &0 &  0 \\
0 &  \mtx{B}_{2,2} &  \mtx{B}_{2,3} & \cdots &0& 0\\
0 & 0 &\mtx{B}_{3,3} & \cdots &0& 0\\
\vdots & \vdots &\vdots &\vdots &\vdots & \vdots \\
0 & 0 & 0 & \cdots & \mtx{B}_{I-1,I} & 0\\
0 & 0 & 0 & \cdots & \mtx{B}_{I,I} &\hat{\mtx{B}}_{I,I+1}\\
\end{bmatrix} \mbox{ with }
\hat{\mtx{B}}_{1,1}= \begin{bmatrix}\mtx{B}_{1,1} & 0\end{bmatrix}
\mbox{ and }
\hat{\mtx{B}}_{I,I+1}=\begin{bmatrix}\mtx{B}_{I,I+1} & 0\end{bmatrix},
$$
where $\mtx{B}_{i,i}$ and $\mtx{B}_{i,i+1}$ correspond to the interaction of the $i^{\rm th}$ interface with the unknowns defined on the proxy circle for the $i^{\rm th}$ and $(i+1)^{\rm th}$ layer respectively.  

The matrix $\hat{\mtx{C}}$ is the matrix  $\begin{bmatrix} 
\mtx{C}\\
\mtx{Z}
\end{bmatrix}$ with its rows reordered as follows
$$
\hat{\mtx{C}} =
\begin{bmatrix}
\hat{\mtx{C}}_{1,1} & 0 & 0 & \cdots & 0 &0\\
\mtx{C}_{2,1} & \mtx{C}_{2,2} & 0 & \cdots &0 &0\\
0 &\mtx{C}_{3,2} & \mtx{C}_{3,3} & \cdots &0 &0\\
\vdots & \vdots &\vdots &\vdots &\vdots & \vdots \\
 0 & 0 & 0 & \cdots & \mtx{C}_{I,I-1} & \mtx{C}_{I,I}\\
 0 & 0 & 0 & \cdots & 0 & \hat{\mtx{C}}_{I+1,I}\\
\end{bmatrix} \mbox{ with }
\hat{\mtx{C}}_{1,1}=\begin{bmatrix}\mtx{C}_{1,1} \\ \mtx{Z}_U\end{bmatrix}
\mbox{ and }
\hat{\mtx{C}}_{I+1,I}=\begin{bmatrix}\mtx{C}_{I+1,I} \\ \mtx{Z}_D\end{bmatrix},
$$
where $\mtx{Z}_U$ is the evaluation of the potential  from the boundary charge density on the first interface on the top of the unit cell, $\mtx{Z}_D$ is the evaluation of the potential from the boundary charge density on the last interface on the bottom of the unit cell,
and $\mtx{C}_{i,i-1}$ and $\mtx{C}_{i,i}$ evaluate the potential due to the charge boundary density on the $(i-1)^{\rm th}$ and $i^{\rm th}$ interface on the
left $L$ and right $R$ vertical walls of the unit cell  in the $i^{\rm th}$ layer to enforce periodicity.

The matrix 
$\hat{\mtx{Q}}$ is a block diagonal matrix obtained from reordering the rows and columns of $ \begin{bmatrix}
\mtx{Q} & \mtx{0}\\
\mtx{V} & \mtx{W}
\end{bmatrix}$.  The block entries are defined as follows
$$
\hat{\mtx{Q}}_{1,1}=\begin{bmatrix}
\mtx{Q}_{1,1} & 0 \\
\mtx{V}_U & \mtx{W}_U
\end{bmatrix}
\hat{\mtx{Q}}_{I+1,I+1}=\begin{bmatrix}
\mtx{Q}_{I+1,I+1} & 0 \\
\mtx{V}_D & \mtx{W}_D
\end{bmatrix}
\mbox{ and }
\hat{\mtx{Q}}_{j,j} = \mtx{Q}_{j,j} \ \forall j = 2,\ldots I.
$$
where $\mtx{V}_U$ and $\mtx{V}_D$ are the only non-trivial blocks in $\mtx{V}$ corresponding to evaluating the potential at $y=y_U$ and $y=y_D$ due to the proxy circle of the top and bottom layer;
$\mtx{W}_U$ and $\mtx{W}_D$ are the non-trivial blocks in $\mtx{W}$ corresponding to evaluating the Rayleigh-Block expansion at $y=y_U$ and $y=y_D$;
$\mtx{Q}_{i,i}$ evaluates the potential on the vertical walls due to the proxy circle of the $i$th layer.

Thanks to the block sparsity pattern of $\mtx{A}$, $\hat{\mtx{B}}$, $\hat{\mtx{C}}$ and $\hat{\mtx{Q}}$ 
the Schur complement matrix $\mtx{S}$ defined 
in (\ref{eq:schur1}) is block tri-diagonal.  The non-zero blocks of $\mtx{S}$ are defined as follows
\begin{equation}
\label{equ:block_def_schur}
\begin{split}
{\mtx{S}}_{1,1}&= \underbrace{\mtx{A}^s_{1,1}+\mtx{A}^{pm}_{1,1}}_{\mtx{A}_{1,1}} - \hat{\mtx{B}}_{1,1}\hat{\mtx{Q}}^\dagger_{1,1} \hat{\mtx{C}}_{1,1}
-\hat{\mtx{B}}_{1,2}\hat{\mtx{Q}}^\dagger_{2,2} \hat{\mtx{C}}_{2,1},\\
{\mtx{S}}_{i,i-1}&= \mtx{A}_{i,i-1} - \mtx{B}_{i,i}\mtx{Q}^\dagger_{i,i} \mtx{C}_{i,i-1},\mbox{ for }i=2,\dots, I\\
{\mtx{S}}_{i,i}&= \underbrace{\mtx{A}^s_{i,i}+\mtx{A}^{pm}_{i,i}}_{\mtx{A}_{i,i}} - \mtx{B}_{i,i}\mtx{Q}^\dagger_{i,i} \mtx{C}_{i,i}
-\mtx{B}_{i,i+1}\mtx{Q}^\dagger_{i+1,i+1} \mtx{C}_{i+1,i},\mbox{ for }i=2,\dots, I-1\\
{\mtx{S}}_{i,i+1}&= \mtx{A}_{i,i+1} - \mtx{B}_{i,i+1}\mtx{Q}^\dagger_{i+1,i+1} \mtx{C}_{i+1,i+1},\mbox{ for }i=2,\dots, I-1\\
{\mtx{S}}_{I,I}&= \underbrace{\mtx{A}^s_{I,I} +\mtx{A}^{pm}_{I,I}}_{\mtx{A}_{I,I}} - \mtx{B}_{I,I}\mtx{Q}^\dagger_{I,I} \mtx{C}_{I,I}
-\hat{\mtx{B}}_{I,I+1}\hat{\mtx{Q}}^\dagger_{I+1,I+1} \hat{\mtx{C}}_{I+1,I}.\\
\end{split}
\end{equation}

\section{The fast direct solver}
\label{sec:FDS}

Recall that the block solve (\ref{eq:blocksolve}) requires first solving for $\hat{\vct{\sigma}}$ and then solving 
for $\vct{x}$.  An efficient way of solving for $\vct{x}$ is to exploit the sparsity pattern in the matrices $\hat{\mtx{Q}}^\dagger$ and $\hat{\mtx{C}}$.  
The entries of $\vct{x}$ are given by 
\begin{equation}
\label{eq:postprocessing_coeff}
\begin{split}
\vct{x}_1 &= - \hat{\mtx{Q}}^\dagger_{1,1} \hat{\mtx{C}}_{1,1} 
\hat{\vct{\sigma}}_1,\\
\vct{x}_i &= - \hat{\mtx{Q}}^\dagger_{i,i} 
\begin{bmatrix}\mtx{C}_{i,i-1} & \mtx{C}_{i,i}\end{bmatrix}
\begin{bmatrix}\hat{\vct{\sigma}}_{i-1}\\
\hat{\vct{\sigma}}_{i}
\end{bmatrix},\mbox{ for }i=2,\dots, I,\\
\vct{x}_{I+1} &= - \hat{\mtx{Q}}^\dagger_{I+1,I+1} \hat{\mtx{C}}_{I+1,I} 
\hat{\vct{\sigma}}_I.\\
\end{split}
\end{equation}

Thus the difficulty in the solution technique lies in solving  
\begin{equation}
\label{eq:schurcompsystem}
{\mtx{S}}\vct{\hat{\sigma}}=\vct{f}.
\end{equation} 

The remainder of this section is dedicated to presenting the construction of the fast direct solver
for the matrix $\mtx{S}$.

\subsection{A closer look at $\mtx{S}$}
\label{sec:lookS}
Recall that $\mtx{S}$ is block tri-diagonal where  the blocks are defined by (\ref{equ:block_def_schur}).  We choose to 
write $\mtx{S}$ as the sum of a full rank block diagonal matrix and a block tri-diagonal matrix whose blocks are low-rank.
Specifically, we express $\mtx{S}$ as
\begin{equation}
\label{eq:shortschurcomp}
\mtx{S} = \mtx{A}_0 +{\mtx{P}},
\end{equation}
where $\mtx{A}_0$ is defined in (\ref{equ:def_Anot}) and $\mtx{P}$ can be defined blockwise as
\begin{equation*}
\begin{split}
\mtx{P}_{1,1}&= \mtx{A}^{pm}_{1,1} - \hat{\mtx{B}}_{1,1}\hat{\mtx{Q}}^\dagger_{1,1} \hat{\mtx{C}}_{1,1}
-\mtx{B}_{1,2}\mtx{Q}^\dagger_{2,2} \mtx{C}_{2,1},\\
\mtx{P}_{i,i-1}&= \mtx{A}_{i,i-1} - \mtx{B}_{i,i}\mtx{Q}^\dagger_{i,i} \mtx{C}_{i,i-1},\mbox{ for }i=2,\dots, I\\
\mtx{P}_{i,i}&= \mtx{A}^{pm}_{i,i} - \mtx{B}_{i,i}\mtx{Q}^\dagger_{i,i} \mtx{C}_{i,i}
-\mtx{B}_{i,i+1}\mtx{Q}^\dagger_{i+1,i+1} \mtx{C}_{i+1,i},\mbox{ for }i=2,\dots, I-1\\
\mtx{P}_{i,i+1}&= \mtx{A}_{i,i+1} - \mtx{B}_{i,i+1}\mtx{Q}^\dagger_{i+1,i+1} \mtx{C}_{i+1,i+1},\mbox{ for }i=1,\dots, I-1\\
\mtx{P}_{I,I}&= \mtx{A}^{pm}_{I,I} - \mtx{B}_{I,I}\mtx{Q}^\dagger_{I,I} \mtx{C}_{I,I}
-\hat{\mtx{B}}_{I,I+1}\hat{\mtx{Q}}^\dagger_{I+1,I+1} \hat{\mtx{C}}_{I+1,I}.\\
\end{split}
\end{equation*}
Since all the blocks in $\mtx{P}$  correspond to non-self interactions, they are low-rank.  Let $\mtx{P}\approx \mtx{LR}$ 
denote the low-rank factorization of $\mtx{P}$.  
Then an approximate solution to (\ref{eq:schurcompsystem}) can be obtained via a Woodbury formula
\begin{equation}
\label{eq:woodbury}
\hat{\bm{\sigma}}=\mtx{S}^{-1}\vct{f}\approx \left( \mtx{A}_0 +\mtx{L}\mtx{R} \right)^{-1}\vct{f} = \mtx{A}_0^{-1}\vct{f} - \mtx{A}_0^{-1}\mtx{L} (\mtx{I}+\mtx{R}\mtx{A}_0^{-1}\mtx{L})^{-1} \mtx{R}\mtx{A}_0^{-1}\vct{f}.
\end{equation} 
%
Since each of the diagonal blocks of $\mtx{A}_0$ is the discretized integral equation on an interface, they are amenable to fast direct solvers
such as the HBS, HSS, HOLDR, etc methods.  Thus $\mtx{A}_0^{-1}$ can be approximated for a cost that scales linearly with respect
to the number of discretization points on the interfaces.  
The remaining complexity lies in creating the low-rank factorization of $\mtx{P}$ and inverting the matrix $\mtx{X}_{\rm wood} = \mtx{I}+\mtx{R}\mtx{A}_0^{-1}\mtx{L}$.

The low-rank factorization of $\mtx{P} $ is handled block wise.  For each non-zero block in $\mtx{P}$, we build the factorization 
for each matrix in the sum independently exploiting the associated physics.  For example, in creating the 
low-rank factorization of the block $\mtx{P}_{i,i}$, we create the low-rank factorization of $\mtx{A}^{pm}_{i,i}$, $\hat{\mtx{B}}_{i,i}\hat{\mtx{Q}}^\dagger_{i,i} \hat{\mtx{C}}_{i,i}$, and $\hat{\mtx{B}}_{i,i+1}\hat{\mtx{Q}}^\dagger_{i+1,i+1} \hat{\mtx{C}}_{i+1,i}$ independently.  Let the low-rank factors be defined as follows: $ \mtx{A}_{i,i}^{pm}\approx \mtx{L}_{i,i}^{pm}\mtx{R}_{i,i}^{pm}$, $\hat{\mtx{B}}_{i,i}\hat{\mtx{Q}}^\dagger_{i,i} \hat{\mtx{C}}_{i,i} \approx \mtx{L}_{i,i,i} \mtx{R}_{i,i,i}$, and $\hat{\mtx{B}}_{i,i+1}\hat{\mtx{Q}}^\dagger_{i+1,i+1} \hat{\mtx{C}}_{i+1,i} \approx \mtx{L}_{i,i+1,i} \mtx{R}_{i,i+1,i}$.
Then $\mtx{P}_{i,i}$ can be approximated by
$$
\mtx{P}_{i,i} \approx 
\begin{bmatrix}
\mtx{L}_{i,i}^{pm} &
\mtx{L}_{i,i,i} & \mtx{L}_{i,i+1,i} 
\end{bmatrix} 
\begin{bmatrix}
\mtx{R}_{i,i}^{pm} \\
-\mtx{R}_{i,i,i} \\ -\mtx{R}_{i,i+1,i} 
\end{bmatrix}.
$$

The technique for creating the low-rank factorization of the matrices $\mtx{A}^{pm}_{i,i}$, $\mtx{A}_{i,i-1}$, and $\mtx{A}_{i,i+1}$ 
is presented in Section 3.1.1 of \cite{multilayer_part1}.  Section \ref{sec:lowrank} presents the technique for creating the 
low-rank factorization of $\hat{\mtx{B}}_{i,i}\hat{\mtx{Q}}^\dagger_{i,i} \hat{\mtx{C}}_{i,i}$.   The low-rank factorization 
of $\hat{\mtx{B}}_{i,i+1}\hat{\mtx{Q}}^\dagger_{i+1,i+1}\hat{\mtx{C}}_{i+1,i}$ can be created in a similar manner.  

The only thing remaining is a fast inversion technique for the matrix $\mtx{X}_{\rm wood}$.  It happens to be the case that this matrix is 
block tri-diagonal. Thus it can be inverted rapidly via a block version of the Thomas algorithm presented in the Appendix of \cite{multilayer_part1}.



\subsection{Low-rank factorization of $\mtx{B}_{i,i}\mtx{Q}^\dagger_{i,i} \mtx{C}_{i,i}$}
\label{sec:lowrank}
Creating the low-rank factorization of $\mtx{B}_{i,i}\mtx{Q}^\dagger_{i,i} \mtx{C}_{i,i}$ requires
dealing with $\mtx{Q}_{i,i}^\dagger$.  
As mentioned in \cite{Cho:15}, filling the entries of the matrix $\mtx{Q}_{i,i}^\dagger$ is not numerically stable. 
To avoid this, we use the truncated SVD-based pseudoinverse.  

\begin{definition}
Let $\mtx{U}_i\mtx{\Sigma}_i\mtx{T}_i^*$ be the SVD of the matrix  $\mtx{Q}_{i,i}$ of size $2M_w\times P$ where 
$\mtx{\Sigma}_i$ is a diagonal rectangular matrix with entries of the singular values of $\mtx{Q}_{i,i}$ and
matrices $\mtx{U}_i$ and $\mtx{T}_i$ are unitary matrices of size $2M_w\times 2M_w$ and $P\times P$, respectively.  
Then the $\epsilon_{\rm Schur}$
-truncated SVD
 is 
$$\mtx{\hat{U}}_i\mtx{\hat{\Sigma}}_i \mtx{\hat{T}}^*_i$$
where $\mtx{\hat{\Sigma}}_i$ is a diagonal square matrix of size $l_i\times l_i$ where $l_i$ is the
number of singular values of $\mtx{Q}_{i,i}$ that are larger than $\epsilon_{\rm Schur}$, $\mtx{\hat{U}}_i$ is a submatrix of size
$2M_w\times l$ of $\mtx{U}_i$ and $\mtx{\hat{T}}_i$ is a submatrix of size $P\times l$ of $\mtx{V}_i$.
\end{definition}

The matrices of the form $\mtx{Q}^\dagger_{i,i} \mtx{C}_{i,i}$ can be approximated by
\begin{equation}
\mtx{Q}^\dagger_{i,i} \mtx{C}_{i,i}
\approx \mtx{\hat{T}}_i\mtx{\hat{\Sigma}}_i^{-1}\mtx{\hat{U}}_i^* \mtx{C}_{i,i}.
\end{equation}

One way of building a low-rank approximation for the matrix product $\mtx{B}_{i,i}\mtx{Q}^\dagger_{i,i} \mtx{C}_{i,i}$ is to let 
\begin{equation}
\label{eq:simple_factors}
\mtx{L}_{i,i,i} = \mtx{B}_{i,i}\mtx{\hat{T}}_i
\mbox{ and }
\mtx{R}_{i,i,i} = \mtx{\hat{\Sigma}}_i^{-1}\mtx{\hat{U}}_i^* \mtx{C}_{i,i}.
\end{equation}
where $\mtx{L}_{i,i,i}$ is of size $N_i\times l_i$, $\mtx{R}_{i,i,i}$ is of size $l_i\times N_i$ and $l_i$ is the number of singular values of $\mtx{Q}_{i,i}$ that are greater than $\epsilon_{Schur}$.  Unfortunately, $l_i$ is far from the optimal rank and thus the resulting 
low-rank factorization of $\mtx{P}_{i,i}$ is far larger than it needs to be.  This has many implications including artificially limiting the 
number of layers that can be simulated on a machine.


To create a closer to optimal rank factorization of $\mtx{B}_{i,i}\mtx{Q}^\dagger_{i,i} \mtx{C}_{i,i}$, we create a low
rank factorization of $\mtx{B}_{i,i}$ and express the low-rank factorization of $\mtx{B}_{i,i}\mtx{Q}^\dagger_{i,i} \mtx{C}_{i,i}$
in terms of those factors.  Let $\mtx{B}_{i,i} \approx \mtx{L}_{B,ii}\mtx{R}_{B,ii}$ denote the low-rank factorization 
of $\mtx{B}_{i,i}$.  Then
\begin{equation}
\mtx{B}_{i,i}\mtx{Q}^\dagger_{i,i} \mtx{C}_{i,i}\approx
\mtx{L}_{B,ii}\mtx{R}_{i,i,i}
\mbox{, with }
\mtx{R}_{i,i,i} = \mtx{R}_{B,ii}\mtx{\hat{T}}_i\mtx{\hat{\Sigma}}_i^{-1}\mtx{\hat{U}}_i^* \mtx{C}_{i,i}.
\end{equation}


In practice, the inner dimensions of this factorization are smaller than $l_i$ as the matrix $\mtx{B}_{i,i}$ 
corresponds to interactions between distant points while the matrix $\mtx{Q}$ is nearly full rank. 
While the second approach requires an additional low-rank factorization, the cost of doing this is more than made up for
by the reduced rank and the fact that these factors can be used for all Bloch phase $\alpha$.  
Additionally, a fast algorithm is used to create the low-rank factorizations of the $\mtx{B}_{i,i}$ and $\mtx{B}_{i,i+1}$ matrices.
Details for constructing that factorization are provided in \cite{zhang2022fast}.

%

%


\section{Extensions}
\label{sec:extensions}

Most applications involve solving (\ref{eq:basic}) for many incident angles.   Additionally, it is often of interest to solve (\ref{eq:basic}) but change the wave number in layer or an interface.  Thus in order for the direct solver to be useful for these applications, it is necessary for it to be able to handle the changes with as minimal work as possible.  Section \ref{sec:bloch} reports how almost all the precomputation can be re-used for solving problems involving many incident angles (which means many different Bloch phases).  Then Section \ref{sec:inter} reports on 
how the solver can be utilized for problems involving changes in a layer whether it is a wave number or an interface geometry.  In both situations, the fast direct solution technique presented in this manuscript scales optimally with the number of discretization points and number of layers.



\subsection{Bloch phase and incident angle dependence}
\label{sec:bloch}
Solving a multilayered media scattering problem for many incidents angles happens frequently in applications.  
For example, in creating a Bragg diagram, the solution for a large collection of incident angles in the range of $-\pi$ to $0$ is needed\cite{2010_chirped}. For each incident angle, there is a corresponding Bloch phase.  Some incident angles share a Bloch phase.  This means that entries in (\ref{eq:rearrangedblock}) change.  Since
in most of these matrices the dependence on Bloch phase is 
a scalar multiplication, the bulk of the precomputation can be reused for all Bloch phases.  For example, this is the case in the matrices that make up $\mtx{A}_{i,i}^{pm}$.  This means that
the low-rank factorization of this matrix can be used for all Bloch phases.  It just needs to be scaled by the Bloch phase $\alpha$ correctly.
An example of a matrix that is incident angle dependent and not just Bloch phase is the matrix $\mtx{W}$.  Thanks to the properties of phase shifts, it is possible to build
$\mtx{W}$ for all incident angles that share a Bloch phase. Section 3.2 of \cite{multilayer_part1} details this incident angle blocking.  
This allows for one direct solver to be built for all incident angles that share a Bloch phase.
A detailed list classifying operations by dependence on Bloch phase and incident angle is given at the beginning of section \ref{sec:numerics}.


\subsection{Changing a layer}
\label{sec:inter}
For applications where there is a desire to solve problems with a change in an interface and/or wave number, the solution technique
presented in this manuscript only requires updating the matrices associated with that layer or interface.  In fact, the cost of updating the 
solver scales linearly with the number of discretization points on the affected interfaces.  This means that updating the solver
is much more efficient than building a new one from scratch.  
Specifically, the speed-up of updating the solver over building a new solver from scratch is the total number of interfaces divided by the number of interfaces changed.  
If for each modified geometry the solution is desired for the same collection of incident angles, all the Bloch phase dependent precomputations that are not related to the changed interfaces or layers can be reused.

\section{Numerical examples}
\label{sec:numerics}
This section illustrates the performance of the proposed fast direct solver for a collection of 
multilayered media scattering problems. Like the solver in \cite{multilayer_part1},
the computational cost of the direct solution technique is broken into four parts; they are:\\
\noindent
\textbf{Precomputation I:}  This consists of all computations for the \textit{fast linear algebra} that are independent of the Bloch phase: the fast application of $\mtx{A}_0^{-1}$, the 
low-rank factors for approximating $\mtx{A}_{i,i}^{pm}$ and $\mtx{A}_{i,j}$ up to scaling by the Bloch phase, 
and the low-rank approximation for the blocks in $\hat{\mtx{B}}$.
The computational cost of this step is $O(N_{total})$ where $N_{total} = \sum_{l=1}^I N_l$, and $N_l$ denotes
the number of discertization points on interface $l$. 

\noindent
\textbf{Precomputation II:} This consists of the remainder of the precomputation that is independent of Bloch phase.
 This includes the evaluation of the Bloch phase-independent components of matrix blocks in $\hat{\mtx{C}}$ and $\hat{\mtx{Q}}$. The computational cost of this step is $O(N_{total})$.
 
  \noindent
  \textbf{Precomputation III:}  This consists of all the precomputation that can be used for incident angles that share a Bloch phase $\alpha$, including scaling matrices by $\alpha$, construction 
of the matrix $\mtx{W}$ accounting for all of the incident angles that share a Bloch phase, constructing the truncated SVD for the diagonal blocks in $\hat{\mtx{Q}}$, combining the low-rank factors for different operators to form the low-rank factors for blocks in the final $\mtx{L}$ and $\mtx{R}$ matrix, and constructing the fast apply
of the Schur complement inverse $\mtx{S}^{-1}$.  Details on the construction of $\mtx{W}$ are provided in section 3.2 of \cite{multilayer_part1}.The computational cost of this step is $O(N_{total})$. 
For a fixed number of discretization points per layer but variable number of layers, the computational cost is $O(I)$.

%
\noindent
\textbf{Solve:} This consists of the application of the precomputed solver for the Schur complement system (\ref{eq:schurcompsystem}) to a right hand side $\vct{f}$ and retrieving the rest of the unknowns via (\ref{eq:postprocessing_coeff}). 
The computational cost of this step is $O(N_{total})$. 
And for a fixed number of discretization points per layer but variable number of layers, the computational cost is $O(I)$.


With the parameters $P_l$, $M_w$, $K$, and $M$ constant,
Precomputation I, II, III and the solve scale linearly with respect to both the number of discretization points per interface and the number of layers (or interfaces). 
This is in contrast to the direct solver from  \cite{multilayer_part1} where Precomputation III and the solve steps have a computational cost that scales cubically with respect to the number of layers (or interfaces).  Thus the solver presented in this manuscript is more efficient for structures with a large number of layers. 

All the results in this section are from implementing the algorithm in MATLAB, except for a Fortran implementation of the interpolatory decomposition used in the HBS compression and low-rank factorizations.   
The experiments were run on a dual 2.3 GHz Intel Xeon Processor E5-2695 v3 desktop workstation with 256 GB of RAM. For all experiments in this section, the parameters $P_l$, $M_w$, $K$, and $M$ are kept fixed.   Specifically,
$P_l = 160$, $M_w = 120$, $K =20$ and $M= 60$.


The experiments considered in this section are the same as the ones considered in \cite{multilayer_part1}.  There are some 
experiments where we were able to consider larger problems thanks to the linear scaling of the new solver with respect to
the number of layers.  For the convenience of the reader, numerical results illustrating the performance of the solver
of \cite{multilayer_part1} are provided in the appendix.

The experiments in section \ref{sec:scaling} illustrate the asymptotic scaling of the proposed fast direct solver.  Section 
\ref{sec:angle_sweep} illustrates the performance of the solver when the solution is desired for many incident angles.
Finally, section \ref{sec:perturbation} illustrates the performance of the solver when a wave number in a layer is changed and when an 
interface is changed.

\subsection{Scaling experiment}
\label{sec:scaling}
\begin{figure}[h]
\centering
\begin{picture}(200,200)(100,0)
\put(50,0){\includegraphics[trim={1 2cm 1 2cm},clip,scale=0.8]{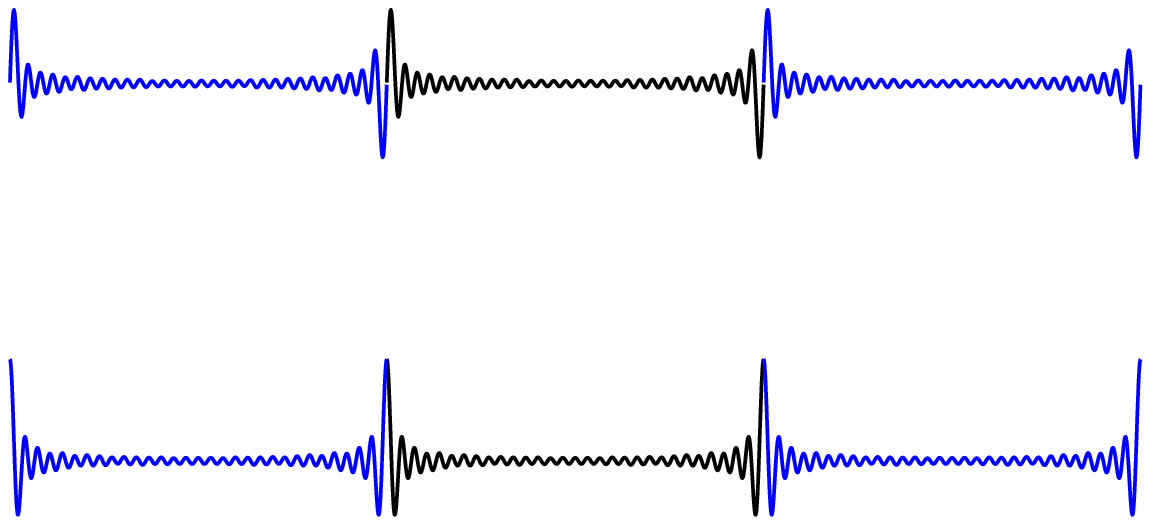}}
\put(210,105){$\gamma_1$ }
\put(210,25){$\gamma_2$ }
\end{picture}
\caption{
Three periods of the interface geometries $\gamma_1$ and $\gamma_2$ as defined in equation (\ref{equ:fouriermap}). 
}
\label{fig:fouriermap_geometry}
\end{figure}

This section illustrates the scaling of the direct solver presented in this manuscript.
We consider a collection of problems where the number of layers 
varies from 3 layers (2 interfaces) to  65 layers (64 interfaces).
The interface geometries are defined by the following two curves $\gamma_1$ and $\gamma_2$ repeated alternatively:  
\begin{equation}\label{equ:fouriermap}
\gamma_1:\,
 \begin{cases}
 x_1(t)=t-0.5\\
 y_1(t) = \frac{1}{60}\sum_{j=1}^{30}  a_j \sin(2\pi j t)\\
 \end{cases}
 \; \mbox{ and } \;
 \gamma_2:\,
 \begin{cases}
 x_2(t)=t-0.5\\
 y_2(t) = \frac{1}{60}\sum_{j=1}^{30}  b_j \cos(2\pi j t)\\
 \end{cases}
\end{equation}
for $t\in[0,1]$,
 where $\{a_j\}_{j=1}^{30}$ and $\{b_j\}_{j=1}^{30}$ are 
 random numbers in $[0,1)$ sorted in descending order.
 Figure \ref{fig:fouriermap_geometry} illustrates the two 
interface geometries.
In each experiment, $\gamma_1$ and $\gamma_2$ are discretized with 
the same number of points $N_i$. The value of $N_i$ is doubled ($1280 \leq N_i\leq 20480$) to demonstrate the linear scaling of the solver.
The wave number of each layer alternates between $10$ and $10\sqrt{2}$. For problems with $64$ interfaces, we exploit the fact that
the interfaces and wave numbers are the same to decrease the memory needed by the solver.  The time in seconds for each part of the precomputation 
using the new solver are presented in Table \ref{tab:timing11}.

\begin{table}[]
 \centering
 \begin{tabular}{|c|c|c|c|c|c|c|}\hline
  &$N_i$ &1280&  2560 & 5120 & 10240 & 20480 \\ \hline \hline
  \multirow{5}{*}{Precomp I} & 2-interface   & 53.4& 111.3&202.6 &350.1 &  594.8 \\
  \cline{2-7}&  4-interface    &101.6 &211.7 & 395.2& 708.2&  1228.8   \\
  \cline{2-7}& 8-interface    &210.5 &431.0 & 793.7&1425.2 &    2457.8   \\
  \cline{2-7}& 16-interface     &442.8 & 898.4&1643.5 &2889.0 &   4914.4   \\
      \cline{2-7}& 32-interface         &861.4& 1678.6&3094.5 & 5725.1&  9489.2*  \\
   \cline{2-7}&  64-interface      & 1721.5* &3611.2* & 6717.0* &12094.9* & 20561.5*    \\
\hline
  \multirow{5}{*}{Precomp II} & 2-interface      &1.7 &2.8 &5.1 &8.6&   16.3      \\
  \cline{2-7}& 4-interface       &2.5 & 4.1& 8.1& 15.1&  29.4       \\
  \cline{2-7}& 8-interface       & 4.2& 8.0&13.9 &26.9 &  49.5     \\
  \cline{2-7}& 16-interface     &7.0 &13.4 &26.1 & 48.2&    116.0   \\
  \cline{2-7}& 32-interface      & 15.2& 26.9&49.1 & 112.2&   217.8*   \\ 
   \cline{2-7}& 64-interface     &34.0*& 57.9*  &107.2* &235.4* & 439.9*     \\
\hline
  \multirow{5}{*}{Precomp III}& 2-interface     & 2.4 & 4.4& 8.9& 15.1&   29.4       \\
  \cline{2-7}& 4-interface       & 4.0&8.1 &17.1 &30.6 &   59.6     \\
  \cline{2-7}& 8-interface      & 7.5& 15.9& 30.3& 59.8&   120.5     \\
  \cline{2-7}& 16-interface     & 14.3 & 34.0& 57.9& 118.4&  252.5    \\
  \cline{2-7}& 32-interface     &27.7  &54.0 &119.1 &252.7 &   505.3 *  \\
   \cline{2-7}& 64-interface     & 52.9* & 121.0* &273.6* &507.8* & 1016.2*     \\
\hline  
  \multirow{5}{*}{Solve} & 2-interface      & 0.1 & 0.3&1.0 &1.6 & 3.8     \\
  \cline{2-7}& 4-interface       &0.3 &0.9 & 1.9& 3.6&  8.4   \\
  \cline{2-7}& 8-interface     & 0.6 &1.6 &4.1 & 9.2&   19.2   \\
  \cline{2-7}&16-interface    & 1.6  &5.3 & 8.8& 19.1&35.6      \\
  \cline{2-7}& 32-interface    & 3.3 &9.2 &19.9 &38.1 &   64.1*  \\ 
   \cline{2-7}& 64-interface    & 8.2*  &15.2* & 32.3*& 64.5*& 137.3*  \\
\hline
  \multirow{5}{*}{Flux error}  & 2-interface       & 1.9e-3 &8.9e-5 & 2.4e-7& 1.1e-9&  1.2e-9      \\
  \cline{2-7}& 4-interface      &3.1e-3 & 1.2e-4 &3.0e-7 &2.7e-9 &    2.6e-9  \\
  \cline{2-7}& 8-interface     & 4.5e-3 &1.8e-4 &1.1e-6 &4.0e-9 &  3.4e-9     \\
  \cline{2-7}& 16-interface     & 8.2e-3 & 2.2e-4& 1.9e-6&7.6e-9 & 6.3e-9     \\
  \cline{2-7}& 32-interface     & 1.7e-2 & 5.2e-4&5.0e-6 & 1.5e-8& 2.0e-8    \\ 
  \cline{2-7}& 64-interface    & 3.0e-2 &8.0e-4 &8.2e-6 &2.8e-8 &  3.4e-8   \\
  \hline
  \end{tabular}
    \caption{Time in seconds and flux error estimates for applying the direct solver to geometries with 2 to 64 interfaces where the interfaces are parameterized smooth curves in Figure \ref{fig:fouriermap_geometry}. The entries marked by ``*" are calculated from a memory efficient implementation of the solver, of which the fast linear algebra for diagonal blocks from repeated interface geometry and wave number is only calculated once and the time is multiplied by the corresponding number of occurrences of the interfaces for the final reported time in the table.
 $N_i$ denotes the number of discretization points for each boundary charge density on the interface.
    The wave number alternates between $10$ and $10\sqrt{2}$.
    }
    \label{tab:timing11}
  \end{table}

For a fixed structure, each part of the 
solution steps scales linearly with respect to $N_i$, the number of discretization points per interface.  
And for fixed $N_i$, each of the steps scales linearly
with respect to the number of interfaces (or the number of layers).
For all tests, Precomputaion I accounts for more than $90\%$ of the total computational cost. Precomputation II accounts for about $5\%$ of the total cost.
Thus the Bloch phase independent parts of the direct solver dominate
the computational cost of building the direct solver.

\begin{remark}
One limitation of the direct solver from \cite{multilayer_part1} is that it cannot handle large number of layers, e.g., 50 layers.
This is because the algorithm in \cite{multilayer_part1} requires taking the pseudoinverse of a matrix whose size scales 
linearly with respect to the number of layers.  This becomes very expensive in terms of memory and number of operations if the number of layers gets large.   The solution technique in this manuscript avoids this problem and only computes pseudoinverses of small block matrices 
defined for each layer thanks to the block diagonal structure of $\hat{\mtx{Q}}$.  
Thus we expect the two solver to have similar performance for a small number of layers but the new solver should be much faster 
for problems with many layers.  

Compared to the scaling results of the original solver in \cite{multilayer_part1} (see Table \ref{tab:old_timing11} in appendix),
when the number of interfaces is less than equal to 8 interfaces, 
the cost per step is similar for all steps except Precomputation III.  The times for Precomputation III using the new solver
are roughly half of the times for the same step with the solver from  
\cite{multilayer_part1}.
This is because the block solve of the new solver does not require applying $\mtx{A}^{-1}$ to $\hat{\mtx{B}}$ when building the Schur complement.

\end{remark}

\subsection{Sweep over multiple incident angles}
\label{sec:angle_sweep}
\begin{figure}[h]
\centering
\begin{picture}(200,250)(100,0)
\put(70,0){\includegraphics[trim={0 0cm 0 0cm},clip,scale=0.5]{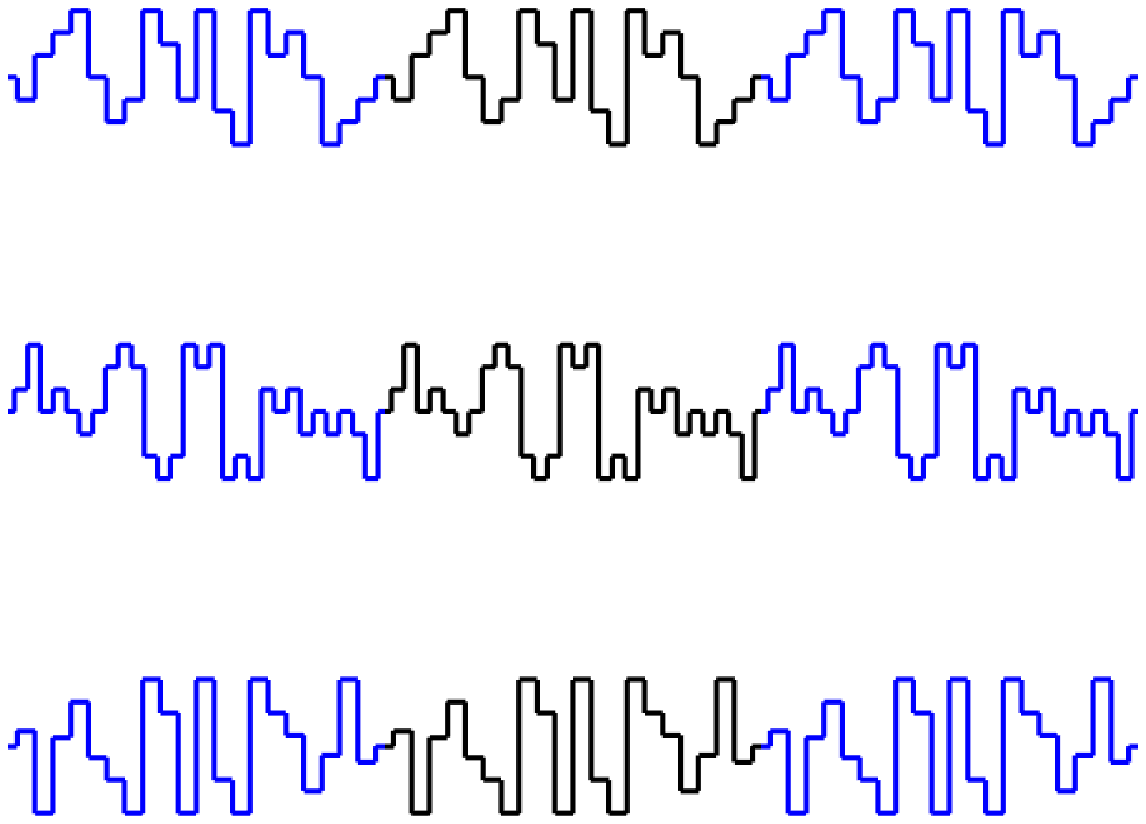}}
\put(70,105){`corner-1' geometry, 42 corners, 10,448 discretization points }
\put(70,55){`corner-2' geometry, 58 corners, 15,104 discretization points }
\put(70,5){`corner-3' geometry, 42 corners, 11,344 discretization points }
\end{picture}
\caption{The three different ``corner'' geometries in the 11-layer structure. 
Three periods are shown. 
See Figure \ref{fig:field_plot} for the full structure. }
\label{fig:corner_geometry}
\end{figure}

\begin{figure}[H]
    \centering
    \includegraphics[scale=.6,trim={1cm 1cm 1cm 1cm},clip]{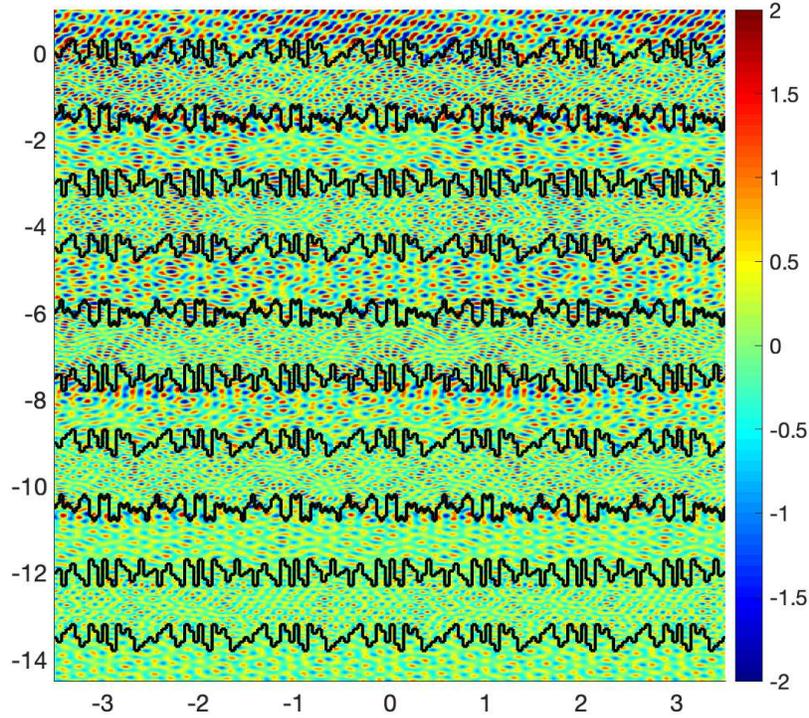}
    \caption{Illustration of the real part of the total field of the solution for the 10-interface structure defined in \cite{multilayer_part1} with incident angle $\theta^{inc}\approx-0.18\pi$, which is at a Wood's anomaly of the top layer. The wave number alternates between $40$ and $40\sqrt{2}$. The flux error is $1.8\times 10^{-8}$.}
    \label{fig:field_plot}
\end{figure}

This section illustrates the performance of the direct solver when the solution 
is desired for many different incident angles.  Specifically, we consider the 
11-layer geometry illustrated in Figure \ref{fig:newinnterface}(a) and approximate the 
solution using the new solver for $287$ different incident angles with $24$ different Bloch phases.
Figure \ref{fig:field_plot} illustrates the real part of the total field for one incident angle.
The interfaces in the 11-layer structure consist of three different
corner geometries, which are referred to as ``corner1", ``corner2" and ``corner3", repeated in order.  
Each of the interfaces contains 40 to 50
right-angle corners. With the five levels of dyadic refinement into
each corner there are 10,000 to 15,000 discretization points per interface.
Figure \ref{fig:corner_geometry} provides more details about the
corner geometries including how many discretization 
points were used on each geometry. 
 The wave number in each layer alternates between $40$ 
and $40\sqrt{2}$.  As stated previously, the solver does group incident angles that share a Bloch phase allowing these incident angles to be solved together.

With this discretization, the average flux error for the $287$ incident angles are $4.5\times10^{-8}$. Table \ref{tab:time_angle} reports the time in seconds for each of the steps.    The incident angle used to create Figure \ref{fig:field_plot} is a Wood's anomaly.  Even in this 
example, the flux error is still on the order of $10^{-8}$, demonstrating that the proposed solution technique is robust at Wood's anomalies.

Table \ref{tab:ols_anglesweep} in the Appendix reports  the performance of the solver from \cite{multilayer_part1} for solving this problem.
As in the previous section, the times for Precomputation I and II for the two solvers are roughly the same.  The
time for Precomputation III for the new solver is twice as fast per Bloch phase as the original solver for this problem.  
This is expected since the new solver no longer requires applying $\mtx{A}^{-1}$ to $\hat{\mtx{B}}$ thus reducing the cost of Precomputation III.  The time for the solve step are essentially the same for both solvers.  

  \begin{table}[]
    \centering
   \begin{tabular}{|c|l|l|l|l|l|}
    \hline
 $N_{total}$& Precomp I & Precomp II & Precomp III & Solve \\
     \hline  
121136 & 2725.4 & 29.8& 2004.5 & 432.2\\
           &            &        &(83.5 per Bloch phase) &( 1.5 per incident angle)\\
\hline
    \end{tabular}
    \caption{Time in seconds for solving 287 incident angles and 24 distinct Bloch phases
    on an 11-layer geometry. The incident angles are sampled from $[-0.89\pi, -0.11\pi]$. 
}
    \label{tab:time_angle}
    \end{table}

\subsection{Local change to the geometry}
\label{sec:perturbation}
\begin{figure}[h]
\centering
\begin{tabular}{cc}
 \includegraphics[trim={5cm .5cm 5cm .5cm},clip,width =4cm]{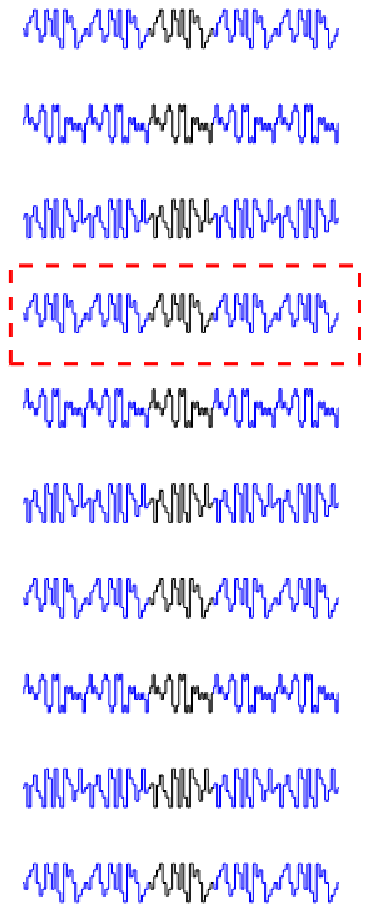}
 &  \includegraphics[trim={5cm .5cm 5cm .5cm},clip,width= 4cm]{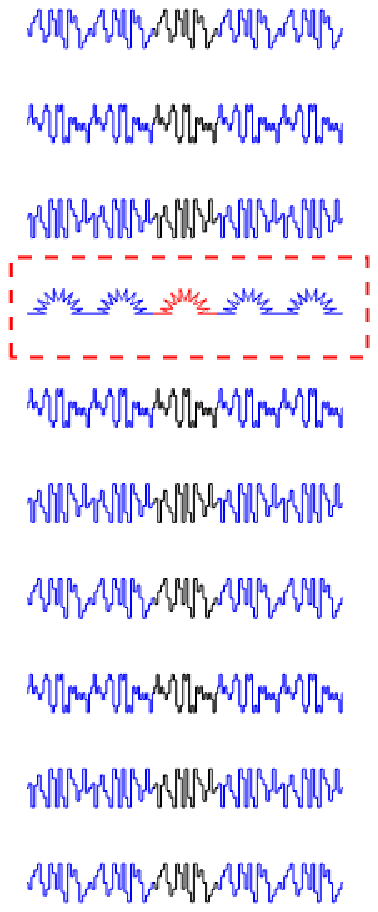}
\end{tabular}
\caption{Illustration of 5 periods of (a) the original 11-layer structure and 
(b) the new structure obtained from replacing the fourth interface with a different geometry. 
The modified interface is in red box. 
}
\label{fig:newinnterface}
\end{figure}

This section illustrates the performance of the proposed direct solver for problems where there is a change in the 
geometry such as changing an interface or the wave number in a layer.  
Specifically, we consider the same 11-layer geometry as in Section \ref{sec:bloch} but change one interface geometry or change the wave number in a layer.
For the change of interface geometry experiment, the fourth interface from the top $\Gamma_4$ is replaced by the ``hedgehog'' geometry as illustrated in Figure \ref{fig:newinnterface}.
The hedgehog geometry consists of 17 sharp corners and cannot
be written as the graph of a function defined on the $x$-axis.  The number of discretization
points on the new interface needed to maintain the same accuracy as the original problem is
 $N_4=14,496$.  
For the change of wave number experiment, the wave number in the second layer from the top of the 11-layer structure $\omega_2$ is changed from $40\sqrt{2}$ to $30$.
For these experiments, the incident angle is fixed at $\theta^{inc}=-\frac{\pi}{5}$. 
Recall that the solver is able to reuse a large portion of the precomputation in these experiments
and only have to process matrices that involve the ``new" portions of the geometry.
Table \ref{tab:new} reports the time in seconds for each step as well as the flux error
for each experiment.  Each step in the precomputation is substantially less expensive.  The smallest
decrease is in Precomputation III which is only a factor $2-3$ faster than the same step when building a new solver
from scratch.  Precomputation I is slower for replacing a wave number than it is for replacing an interface because 
replacing the wave number involves changing matrices for two interfaces.  Even so, building the solver using the 
method presented in this manuscript for the case of replacing a wave number is roughly $5.3$ times faster than building a 
new solver from scratch for the problem.  There is roughly a $8.9$ times speed-up for using the solver presented in this manuscript for the replaced interface problem
instead of building a solver from scratch.  These speed-up numbers will be even greater for more than one Bloch phase. 
In comparison with the solver in \cite{multilayer_part1} for these problems (Table \ref{tab:old_perturb} in the Appendix), the new solver 
is faster in Precomputation III.  

\begin{table}[h]
\centering
\begin{tabular}{|c|c|c|c|
}
\hline
&  Original problem &Replace interface $\Gamma_4$ & Change wave number $\omega_2= 30$ \\ \hline
$N_{total}$&121,136 &   125,184  &121,136 \\ \hline
Precomp I &2305.8& 232.3 &454.9\\ \hline
Precomp II &37.4&8.3  & 3.9\\ \hline
Precomp III  &86.1& 30.2 &  27.1 \\ \hline
Solve &15.0&7.3& 6.9 \\ 
\hline
Flus error &1.3e-8 &3.0e-8 &5.4e-8\\
\hline 
\end{tabular}
\caption{Time in seconds for constructing and applying the fast direct solver 
to an 11-layer geometry (first column), a geometry that 
has the fourth interface changed (second column) and the wave number for the second layer 
changed from $40\sqrt{2}$ to $30$ (third column). $N_{total}$ is the number of discretization points on the 
interfaces in the unit cell.  }
\label{tab:new}
\end{table}%

\section{Conclusion}
\label{sec:summary}
This manuscript presents a new fast direct solver for two-dimensional quasi-periodic scattering problems in multilayered structures. The solution technique is based upon the 
robust integral equation formulation from \cite{Cho:15}.  The solver in this paper approaches the block solve of the linear system in the same manner as in
 \cite{Cho:15} with the addition of physics based accelerated linear algebra.  The resulting solver is more efficient than the solver proposed in \cite{multilayer_part1} for problems with many layers.  For low frequency problems, the solver presented in this manuscript 
 scales linearly with respect to the number of discretization points per interface as well as the number of interfaces. Thus it will be useful in applications
 involving complicated interfaces and many layers.  Numerical results illustrated that for a geometry with eleven layers requiring over 200 solves, the fast direct solver in this manuscript is about 135 times faster than building a new solver from scratch for each right hand side.  For problems where there are changes in a layer, whether it is wave number or the interface geometry, the 
 proposed solver can be updated by a cost that scales linearly with respect to the number of discretization points on the interfaces affected by the update.  For future directions, we are interested in extending the work to higher frequency and three dimensions.

\appendix
\section{Numerical Results for the algorithm presented in \cite{multilayer_part1}}
For fair comparison, a new implementation of the solver in \cite{multilayer_part1} is rerun on the same desktop used for the numerical section of this manuscript. The results are summarized below.

The direct solver presented in \cite{multilayer_part1} 
solves the block system (\ref{eq:bigblock}) via the 
following block solve
\begin{equation}
    \begin{split}
    \vct{x}&=-\left(\hat{\mtx{Q}}-\hat{\mtx{C}}\mtx{A}^{-1}\hat{\mtx{B}}\right)^\dagger\hat{\mtx{C}}\mtx{A}^{-1}\vct{f}\\
    \hat{\vct{\sigma}}&=-\mtx{A}^{-1}\hat{\mtx{B}}\vct{x}+\mtx{A}^{-1}\vct{f}
    \end{split}
\end{equation}
Note that this is a different processing of the system than 
what was presented in this manuscript.  
The computational cost of the solution technique is reported in four categories:
\begin{itemize}
\item Precomputation I:  This consists of all computations for the fast linear algebra that are independent of Bloch phase. This includes the fast application of $\mtx{A}_0^{-1}$, and the 
low rank factors $\mtx{L}_{ij}$ and $\mtx{R}_{ij}$ needed to make $\mtx{L}$ and $\mtx{R}$. 
 \item Precomputation II: This consists of the remainder of the precomputation that is independent of Bloch phase. 
\item Precomputation III:  This consists of all the precomputation that can be used for incident angles that share a Bloch phase $\alpha$, including scaling matrices by $\alpha$, construction 
of the matrix $\mtx{W}$, constructing the fast apply
of $\mtx{A}^{-1}$, evaluating the Schur complement matrix $\mtx{S}=-\left(\hat{\mtx{Q}}-\hat{\mtx{C}}\mtx{A}^{-1}\hat{\mtx{B}}\right)$ , and constructing 
the pseudoinverse of the Schur complement matrix f $\mtx{S}$ via $\epsilon_{\rm Schur}$ SVD . 
\item Solve: This consists of the application of the precomputed solver to the right hand side to evaluate  $\bm{x}$ and then $\hat{\bm{\sigma}}$. 
\end{itemize}
Let ${N}_{total}$ be the total number of discretization points for all interfaces  and $I$ be the number of interfaces. For a fixed structure (with fixed number of interfaces), the cost of Precomputation I, II, III and the solve is $O({N}_{total})$.
For a structure with variable number of interfaces but each interface is discretized with a fixed number of points, the cost of Precomputation I and II is $O(I)$ while that of Precomputation III and the solve is $O(I^3)$.


The results for the scaling experiment are given in Table \ref{tab:old_timing11}. Table \ref{tab:ols_anglesweep} reports the angle sweeping tests, and Table \ref{tab:old_perturb} illustrates the results for modifying the solver to accommodate an interface change and a layer wave number change.  The tested structures are designed to be the same as the ones described in section \ref{sec:scaling}, \ref{sec:angle_sweep} and \ref{sec:perturbation} or the original numerical section of \cite{multilayer_part1}.
\begin{table}[H]
 \centering
 \begin{tabular}{|c|c|c|c|c|c|c|}\hline
  &$N_i$ & 1280& 2560 & 5120 & 10240 & 20480 \\ \hline \hline
  \multirow{3}{*}{Precomp I} & 2-interface &45.4 &92.1 &172.0 &319.6 &   547.8 \\
  \cline{2-7}& 4-interface  & 94.9&195.3 &366.9 & 667.8& 1151.1 \\
  \cline{2-7}& 8-interface & 198.3& 400.3& 737.6& 1337.0 & 2297.6  \\
  \hline
  \multirow{3}{*}{Precomp II} & 2-interface  & 1.3& 2.1& 4.5&  8.6 & 16.9 \\
  \cline{2-7}& 4-interface  &1.9 & 3.6& 7.2&  14.5 & 29.7  \\
  \cline{2-7}& 8-interface  & 4.0&7.2 & 12.9& 31.1 &60.7 \\\hline
  \multirow{3}{*}{Precomp III}& 2-interface  &2.0 &5.0 &10.8 &  22.9  & 43.4  \\
  \cline{2-7}& 4-interface  &4.9 & 12.6& 23.4&  49.2 &  93.9    \\
  \cline{2-7}& 8-interface &14.5 &30.2 & 58.5& 118.9  & 233.4      \\\hline  
  \multirow{3}{*}{Solve} & 2-interface  & 0.1& 0.5& 1.9& 3.2  &3.3     \\
  \cline{2-7}& 4-interface  &0.7 & 1.3& 3.8& 8.2  & 15.7  \\
  \cline{2-7}& 8-interface &1.9 &4.0 & 8.2& 11.9  & 29.6  \\\hline
  \multirow{3}{*}{Flux error} & 2-interface  &4.2e-5 &6.9e-6 & 2.3e-8&  3.8e-10 &    4.5e-10 \\
  \cline{2-7}& 4-interface  &9.8e-5 & 8.0e-6& 8.9e-8& 4.1e-10  & 7.7e-10    \\
  \cline{2-7}& 8-interface &2.1e-4 & 1.2e-5& 1.5e-7  & 4.6e-11 &4.6e-10     \\\hline
  \end{tabular}
    \caption{Time in seconds and flux error estimates for applying the direct solver in \cite{multilayer_part1} to geometries with 2 to 8 interfaces where the interfaces are parameterized smooth curves  defined in section 4.1 of \cite{multilayer_part1}.
 $N_i$ denotes the number of discretization points for each boundary charge density on the interface.
    The wave number alternates between $10$ and $10\sqrt{2}$.
    }
    \label{tab:old_timing11}
  \end{table}
  
    \begin{table}[H]
    \centering
   \begin{tabular}{|c|l|l|l|l|l|}
    \hline
 $N_{total}$& Precomp I & Precomp II & Precomp III & Solve \\
     \hline  
121136 & 2490.0 & 31.0& 4323.6 & 456.7\\
           &            &        &(180.2 per Bloch phase) &( 1.6 per incident angle)\\
\hline
    \end{tabular}
    \caption{Time in seconds for solving 287 incident angles and 24 distinct Bloch phases
    on an 11-layer geometry via the direct solver given in \cite{multilayer_part1}. The incident angles are sampled from $[-0.89\pi, -0.11\pi]$. The average flux error for all solved incident angles is 1.7e-8.
}
    \label{tab:ols_anglesweep}
    \end{table}

  \begin{table}[H]
\centering
\begin{tabular}{|c|c|c|c|
}
\hline
&  Original problem &Replace interface $\Gamma_4$ & Change wave number $\omega_2= 30$ \\ \hline
$N_{total}$&121,136 &   125,184  &121,136 \\ \hline
Precomp I &2320.0& 226.2 & 440.2 \\ \hline
Precomp II &37.0&8.4  & 4.0\\ \hline
Precomp III  &110.1& 30.2 & 107.2  \\ \hline
Solve &19.1&12.7& 11.4 \\ 
\hline
Flux error & 3.4e-8 & 4.2e-9 &4.0e-9\\
\hline 
\end{tabular}
\caption{Time in seconds for constructing and applying the fast direct solver in \cite{multilayer_part1} 
to an 11-layer geometry (first column), a geometry that 
has the fourth interface changed (second column) and the wave number for the second layer 
changed from $40\sqrt{2}$ to $30$ (third column). $N_{total}$ is the number of discretization points on the 
interfaces in the unit cell.  }
\label{tab:old_perturb}
\end{table}%
\end{document}